\documentclass[11pt]{article}

\usepackage{amssymb}
\usepackage{graphicx}
\usepackage{xcolor}
\usepackage{colortbl}

\newtheorem{theorem}{Theorem}
\newtheorem{lemma}{Lemma}[theorem]

\newcommand{\abs}[1]{|{#1}|}

\begin{document}
\bibliographystyle{plain}

\title{ %Estimating the Smooth Counting Function in Practice}
	An Algorithm for Ennola's Second Theorem 
  and Counting Smooth Numbers in Practice} 

\author{Chloe Makdad \\
Computer Science \&\ Software Engineering \\
Mathematics, Actuarial Science, \&\ Statistics \\
  Butler University, Indianapolis IN, USA \\
\texttt{cmakdad@butler.edu}
\and
Jonathan P.~Sorenson \\
Computer Science \&\ Software Engineering \\
  Butler University, Indianapolis IN, USA \\
  \texttt{jsorenso@butler.edu} }
\date{\today}

\maketitle

\begin{abstract}
Let $\Psi(x,y)$ count the number of positive integers $n\le x$
  such that every prime divisor of $n$ is at most $y$.
There are a number of applications where values of $\Psi(x,y)$
  are needed, such as in optimizing integer factoring and discrete
	logarithm algorithms \cite{CP,MOV} and generating factored smooth
	numbers uniformly at random \cite{BS2020}. 
Note that such numbers % uniformly generated factored smooth numbers 
  are useful in at least one post-quantum cryptography protocol
  \cite{Couveignes2006,RS2006}.

Given inputs $x$ and $y$, what is the best way to estimate $\Psi(x,y)$?
We address this problem in three ways:
  with a new algorithm to estimate $\Psi(x,y)$,
  with a performance improvement to an established algorithm,
  and with empirically based advice on how to choose an algorithm
  to estimate $\Psi$ for the given inputs.

Our new algorithm to estimate $\Psi(x,y)$
  is based on Ennola's second theorem \cite{Ennola69},
  which applies when $y< (\log x)^{3/4-\epsilon}$ for $\epsilon>0$.
It takes $O(y^2/\log y)$ arithmetic operations of precomputation
  and $O(y\log y)$ operations per evaluation of $\Psi$.

We show how to speed up Algorithm HT \cite{HS97}, which is based on
  the saddle-point method of Hildebrand and Tenenbaum \cite{HT86},
  by a factor proportional to $\log\log x$, by applying Newton's method
  in a new way.

And finally we give our empirical advice based on five algorithms
  to compute estimates for $\Psi(x,y)$.
The challenge here is that the boundaries of the ranges of applicability,
  as given in theorems, often include unknown constants or small values
  of $\epsilon>0$, for example, that cannot be programmed directly.
\end{abstract}

%%%intro

\section{Introduction\label{sec:intro}}

\nocite{CP,BS,MOV,Granville2008,HT93,HT86,Tenenbaum,RS2006,Couveignes2006}

Let $\Psi(x,y)$ count the number of integers $n\le x$ such that the
  largest prime divisor of $n$ is $\le y$.
There are a variety of algorithms to estimate the value of $\Psi(x,y)$
  in the literature
  \cite{BS13,Bernstein95-2,Bernstein98,Bernstein2002,
  Ennola69,HS97,KTP76,PS06,
  Sorenson2000,Suzuki2004,Suzuki2006,vdLW69}.
All these methods have one drawback or another.
Some have very slow runtimes,
  some are inaccurate in practice, and some have very limited
  ranges of applicability.
To make this even more difficult,
  in most cases, the theorems that provide a
  region of the $x/y$ plain where the algorithm's accuracy has a guarantee
  is not specified explicitly.
This makes the boundaries of such regions impossible to program.

In this paper, we try to determine the best way to estimate $\Psi(x,y)$
  for specific values of $x$ and $y$ in practice.
To do this, we implemented many algorithms for $\Psi(x,y)$ and
  explored their boundaries of applicability empirically.
We present these results in \S\ref{sec:eval}.

In the process of our study, we noticed that the second theorem
  in Ennola's paper \cite{Ennola69} had not been tried, to our knowledge.
The first theorem in that paper is well known and is quoted,
  for example, in \cite[\S5.2]{Tenenbaum}.
We found that the range of applicability of this second theorem, in practice,
  far, far exceeds its proven guarantee.
So we begin with an exposition of our algorithm based on Ennola's second
  theorem, and we analyze its running time below in \S\ref{sec:ennola}.
We believe this algorithm is completely new.

We also show, in \S\ref{sec:htalpha}, how to trim a factor proportional
  to $\log \log x$ from the running time of Algorithm HT \cite{HS97}.

We conclude in \S\ref{sec:wrapup} with some comments.

%%%%ennola

\section{An Algorithm Based on Ennola's Second Theorem\label{sec:ennola}}

We begin this section by reviewing Ennola's second theorem, 
  then we present the algorithm,
  we give an analysis of the running time and the space used,
  and we conclude with some practical notes and data on
  the algorithm's accuracy in practice.

\subsection{Ennola's Theorem}

%let $\Psi(x,y)$ be the number of integers $\le x$ that are composed entirely
%  of primes $\le y$, that is, $y$-smooth.
%There are a number of estimates for $\Psi(x,y)$ in the literature;

Let $2\le y <x$.
Ennola's theorem applies when $y$ is very small, namely when 
  $y\le (\log x)^{3/4-\epsilon}$ for $\epsilon>0$.

We need a few definitions.

First, define the sequence $c_n$ with $c_0=1$, $c_1=1/2$, $c_{2k+1}=0$ for
  $k\ge 1$, and
\begin{equation}
 c_{2k} = (-1)^{k+1} \frac{2\zeta(2k)}{(2\pi)^{2k}} \label{eq:ck} 
\end{equation}
for $k\ge 1$; here $\zeta$ is the Riemann zeta function.

Next, define $d_n$ to be 
  the coefficient of $s^n$ in the following power series:
\begin{equation}
  \prod_{p\le y} \left(\sum_{k=0}^\infty c_k (\log p)^k s^k \right)
  \label{eq:dnpoly} .
\end{equation}
Finally, define
$$
  R(t) := \sum_{n=0}^{\pi(y)} \frac{ d_n t^n}{(\pi(y)-n)!}.
$$

\begin{theorem}[Ennola \cite{Ennola69}]
  Let $2\le y \le (\log x)^\theta$ with $0<\theta<3/4$.
Then
$$
  \Psi(x,y) = \left(\prod_{p\le y} \frac{\log x}{\log p} \right)
    R( 1/\log x )
   \left(1+ O( (\log x)^{-1/8+\epsilon}) \right)
$$
\end{theorem}

Ennola also gave an upper bound for the tail of the sum defining $R(t)$,
  which will enable us to compute with only the first few terms. 
Let $m$ be an integer with $2\le m \le \pi(y)$.
In practice, we ended up using $m=\pi(y)$ almost always.

Define
\begin{equation}
  R_m(t) := \sum_{n=0}^{m} \frac{ d_n t^n}{(\pi(y)-n)!}.  \label{eq:rn}
\end{equation}
Then $R(t)=R_{\pi(y)}(t)$.
Ennola \cite[(27), p.9]{Ennola69} proved
\begin{lemma}
Let $m\ge n_0$ where $n_0=\max\{e\log y, ey^2/(\log x \log y)\}$.
If $0<t<1$, then
\begin{eqnarray*}
  R(t)-R_m(t) &=& 
    \sum_{n=m+1}^{\pi(y)} \frac{ d_n t^n}{(\pi(y)-n)!} \\
  &\le&
    \sum_{n=m+1}^{\pi(y)} \frac{ |d_n| t^n}{(\pi(y)-n)!} \\
    &<& \frac{1}{2^{m}\pi(y)!}
\end{eqnarray*}
\end{lemma}
Note that since $d_0=1$, the first term of the sum is $1/\pi(y)!$,
  and in fact Ennola showed that $R(1/\log x)\gg 1/\pi(y)!$
  \cite[(23),(28)]{Ennola69}.
Given this, it is not too surprising that
  Ennola was able to show that if $m\ge n_0$ then we have
$$
  R(1/\log x) = R_m(1/\log x) (1+O(2^{-m})).
$$

Below, we will choose $m\ge\log_2\log x$,
  allowing us to compute $R_m$ in place of $R$ with no ill effect
  on the overall relative error.

\subsection{The Algorithm}

\subsubsection{Precomputation}

Computing the $d_n$ coefficients is the bottleneck of the algorithm,
but we can precompute them since they are the same for all $(x,y)$
input pairs with the same $y$ value.
And in fact, if we know beforehand the maximum $y$ we will need to
  accept as input, we can find the $d_n$ coefficients for all $y^\prime\le y$
  along the way for no extra cost,
  aside from the space needed to store the coefficients.

Normally we use $m=\pi(y)$ for precomputation, unless we have a
  reasonably tight range of values for $x$ that we can bound beforehand.
We address this situation further below and give a tighter runtime
  analysis for precomputation for when this is the case.

\begin{enumerate}
\item Compute the list of primes up to $y$ using a sieve.
%\item Set 
%  $m:=\lfloor
%     \min\{ \pi(y), \max\{\log_2\log x, e\log y, ey^2/(\log x\log y) \}\}
%       \rfloor$.
\item Compute $\zeta(2k)$ for integers $k$, $0<2k\le m$.
  We do this by noting that $\zeta(2)=\pi^2/6$ %,\zeta(4)=\pi^4/90$,
  and using the recurrence
  $$
    \zeta(2k)=\frac{1}{k+1/2} \sum_{j=1}^{k-1} \zeta(2j)\zeta(2k-2j).
  $$
\item Compute the $c_n$ for $0\le n\le m$ using (\ref{eq:ck}).
\item Next, we compute the $d_n$ for $0\le n\le m$.
  Recall that $d_k$ is defined as the coefficient of $s^k$ 
    in (\ref{eq:dnpoly}).

  Define the polynomial 
  $f_{p,m}(x):= \sum_{k=0}^m c_k (\log p)^k x^k $ 
  for $p$ a prime.
  Compute $F$ as follows:
  \begin{tabbing}MM\=MM\=MM\=\kill
     \>$F:=f_{2,m}(x)$;\+ \\
     for each prime $p$, $2<p\le y$: \+\\
        $F:=F\cdot f_{p,m}(x) \bmod x^{m+1}$; 
  \end{tabbing}
  The coefficients $d_n$ for $n\le m$ are then the coefficients of $F$.
\end{enumerate}
To save the $d_k$ coefficients for $y^\prime<y$, one simply pulls their
  values off $F$ in the last step immediately after all primes $p\le y^\prime$
  have been processed.

\subsubsection{The Algorithm}
With the $d_n$ coefficients precomputed, the algorithm is as follows.
\begin{enumerate}

\item Set
  $m:=\lfloor
     \min\{ \pi(y), \max\{\log_2\log x, e\log y, ey^2/(\log x\log y) \}\}
       \rfloor$.

\item Compute $R_m(1/\log x)$ using (\ref{eq:rn}).

\item Compute $T:=\prod_{p\le y} \frac{\log x}{\log p}$
   and output the estimate $T\cdot R_m(1/\log x)$.
  
\end{enumerate}

\subsection{Complexity Analysis}

We assume basic arithmetic operations on integers and floating-point numbers,
  such as addition and multiplication, take constant time.
We also assume that basic special functions, like $\log$ and $\exp$
  take constant time.
In practice, we used the standard \texttt{long double} data type in C/C++.

To maintain a relative error of $1+O(1/\log x)$, we need at least 
  $O(\log\log x)$ bits of precision in all our calculations.
We will measure space in the number of machine words, under the assumption that
  each word holds one floating point number of the necessary precision.

\subsubsection{Precomputation}

\begin{enumerate}
  \item The primes up to $y$ can be found using $O(y/\log\log y)$ arithmetic
     operations; see \cite{AB2004}.
     Storing these primes takes $O(y/\log y)$ words of space.
%  \item This takes constant time and space.
  \item This takes $\sum_{k=1}^{m/2} (k-1) = O(m^2)$ time.
     $O(m)$ words are required to store the $\zeta$ function values.
  \item This is $O(m)$ operations if you are careful about how the powers of
     $2\pi$ are computed.
     Again, $O(m)$ words of space are needed for the $c_n$ values.
  \item 
     Computing each $f_{p,m}(x)$ takes $O(m)$ operations, if the terms
       are computed from low degree to high.

     $O(m\log m)$ operations are needed for the convolutions 
     to compute a single polynomial product, using FFT techniques.  

     The total time, then, is $O(\pi(y) m\log m)$ for this step.
     It uses $O(m)$ words of space.

     In practice we used a simple $O(m^2)$ algorithm for polynomial
       multiplication.
\end{enumerate}
The total time for precomputation is $O(\pi(y) m\log m)$ or,
  when $m=\pi(y)$,
  $O( \pi(y)^2 \log y) = O( y^2/\log y)$ for all precomputation up to $y$.

The total space is $O( \pi(y)^2)$ to store all the $d_n$ coefficients for
  all $y^\prime \le y$.

\subsubsection{The Algorithm}
\begin{enumerate}
  \item This is constant time and space.
   \item Computing $R_m(1/\log x)$ takes $O(m+\pi(y))$ operations.
     Note that powers of $1/\log x$ and the factorial denominators should be
     computed in opposite directions first.

     Again, $O(m)$ words of space suffice for this step.

   \item Computing $T$ takes $O(\pi(y))$ operations and constant space.
\end{enumerate}
So after precomputation, the time is $O(\pi(y))$ operations
 to compute an estimate for $\Psi(x,y)$, independent of $x$ (or $m$).

If we know $x$, or have a bound on its range relative to $y$,
  then computing $m$ may save a bit of precomputation time in some cases,
  as the following table shows.
%The time for the polynomial products needed to compute the $d_n$ values 
%  is the bottleneck of this algorithm, 
%  at $O(\pi(y) m\log m)$ total operations.
%, which works out to $O(y^2/\log y)$ in the worst case
%The total space used is $O(m+\pi(y))$ words of space.
%
%The running time of
%  $O(\pi(y) m\log m)$ is, in the worst case, $O(y^2/\log y)$, but can
%  be a bit better than that, depending on how large $y$ is compared to $x$:
  \\
\ \\
\begin{tabular}{l|l|l}
\textbf{Ranges for $x,y$} & \textbf{$m$} & \textbf{Time (ops)} \\ \hline
$\pi(y)\le \log_2\log x$ & $m=\pi(y)$ &
   $O( (\log\log x)^2 \log\log\log x)$ \\\hline
$\log_2\log x < \pi(y) \le \sqrt{\log x}$ &
  $m=\log_2\log x$ &
  $O( \pi(y) \log\log x \log\log\log x)$ or \\
	&($e\log y$ is smaller) & $O(\sqrt{\log x} \log\log x\log\log\log x)$ \\ \hline
$\sqrt{\log x}<\pi(y)\le \pi(\log x)$ &
  $m=ey^2/(\log x \log y)$ & $O(y^3/(\log x \log y))$ \\
%  && \\ 
	\hline
$\log x \ll y$ & $m=\pi(y)$ & $O(y^2/\log y)$ \\
\end{tabular}
\\
\ \\
All cases are bounded by $O(y^2/\log y)$.
Note that the last row of the table is outside the guaranteed
  range given in Ennola's second theorem.

\subsection{Practical Notes}

\subsubsection{Computing $T\cdot R_m(1/\log x)$}

In practice, computing steps 2 and 3 of the main algorithm, especially 2,
  can lead to overflow or underflow if using fixed precision floating point
  numbers, such as the \texttt{long double} datatype in C++.
In addition to this, if one is careful, it is possible to evaluate
  $R$ in time linear in $\pi(y)$, as we will now show.

Define 
$$
  f_n := \frac{ (\log x)^{\pi(y)-n} }{ (\pi(y)-n)! } .
$$
Then we have $f_n = f_{n+1} \cdot (\log x)/(\pi(y)-n)$ when $0\le n<\pi(y)$,
  and $f_{\pi(y)}=1$.
This gives us
\begin{eqnarray*}
  T\cdot R_m(1/\log x) &=&
   \prod_{p\le y} \frac{\log x}{\log p}
   \sum_{n=0}^m \frac{d_n (\log x)^{-n}}{(\pi(y)-n)!} \\
  &=&
   \prod_{p\le y} \frac{1}{\log p}
   \sum_{n=0}^m \frac{d_n (\log x)^{\pi(y)-n}}{(\pi(y)-n)!} \\
  &=&
   \prod_{p\le y} \frac{1}{\log p}
   \sum_{n=0}^m d_n\cdot f_n.
\end{eqnarray*}
The following pseudocode fragment will compute this:
  %$sum := T\cdot R(1/\log x)$:
\begin{tabbing}MM\=MM\=MM\=\kill
\>$fn:=1$; \+\\
  if $m=\pi(y)$ \+\\
    then $sum:=d_{\pi(y)}\cdot fn$\\ 
    else $sum:=0;$ \-\\
  endif; \\
  for $n=\pi(y)-1$ downto $0$ do:\+\\
    $fn:=fn\cdot (\log x)/(\pi(y)-n)$; \\
    if $n\le m$ then $sum:=sum+d_n\cdot fn$; endif;\- \\
  endfor; \\
  $P:=1$; \\
  for each prime $p\le y$ do:\+\\
    $P:=P\cdot 1/(\log p)$; \-\\
  endfor; \\
  output $sum\cdot P$;
\end{tabbing}

\subsubsection*{Accuracy}

We conclude this section with empirical results on the accuracy
  of this new algorithm.
In the table below we give the ratio of the value given by
  our new algorithm over the exact value of $\Psi(x,y)$ for
  various values of $x,y$.
\begin{center}
\begin{tabular}{r|lllll}
  \hline
$y$ & $x=2^{15}$ &$x=2^{20}$ &$x=2^{25}$ &$x=2^{30}$ &$x=2^{33}$ \\
  \hline
   $2^5$ &
0.999972 & 1.00009 & 1.00001 & 0.999978 & 0.999998 \\
$2^{10}$ &
1.0083 & 0.995777 & 1.00115 & 0.999969 & 1.00052 \\
$2^{15}$ &
	-- & 1.00219 & 1.00472 & 0.994501 & 1.00183 \\
  \hline
\end{tabular}
\end{center}
Note that when $y=32=2^5$, for Ennola's second theorem to apply,
  we would require that $y< (\log x)^{3/4}$;
  this would imply $x>\exp 32^{4/3}$, a 44-digit number.
So this method seems to apply to a much wider $x,y$ range than is currently
  proven.
And, although its preprocessing makes it very slow for larger $y$,
  it seems to be as accurate, if not more accurate, than 
  Algorithm HT \cite{HS97}.

%%%%htalpha

\section{Improving Algorithm HT\label{sec:htalpha}}

At a high level, Algorithm HT \cite{HS97} uses Newton's method to find
  the zero, $\alpha$, of a continuous function.
Our idea to improve the algorithm is to first find an approximation
  to $\alpha$, called $\alpha_f$, using the version of Algorithm HT
  that assumes the Riemann Hypothesis to bound the error when
  estimating the distribution of
  primes, allowing for faster summing of functions of primes 
  \cite{Sorenson2000}.
Then, starting from $\alpha_f$, Newton's method is applied in the
  context of the original Algorithm HT, allowing for much faster
  convergence, often requiring only one iteration in practice,
  yet providing the same level of accuracy as Algorithm HT.

We begin this section with a review of Algorithms HT and HT-fast,
then present our new twist, Algorithm HT$\alpha$,
and wrap up the section with some implementation results.

\subsection{Algorithm HT and Algorithm HT-fast}

A theorem from Hildebrand and Tenenbaum \cite{HT86} gives us
    \[\Psi(x,y) \approx HT(x,y,\alpha(x,y))\]
    uniformly for $2\leq y\leq x$, where
    \[HT(x,y,s) := \frac{x^s\zeta(s,y)}{s\sqrt{2\pi\phi_2(s,y)}}\]
    and $\alpha(x,y)$ is the unique solution to
    \[\phi_1(\alpha,y)+\log x = 0.\]
      Here,
      \begin{eqnarray*}
	      \zeta(s,y) &=& \prod_{p\leq y} (1-p^{-s})^{-1}, \\
	      \phi_1(s,y)& =& -\sum_{p\leq y} \frac{\log p}{p^s-1}, \\
	      \phi_2(s,y) &=& \sum_{p\leq y}\frac{p^s(\log p)^2}{(p^s-1)^2}.
      \end{eqnarray*}

{Algorithm $HT$} \cite{HS97}:
\begin{enumerate}
    \item Find all primes $p\leq y$.
    \item Starting at $\alpha_0 := \log(1+y/(5\log x))/\log y$, 
	   approximate $\alpha$ by $\alpha'$, 
	  where $s=\alpha$ is the solution to $f(s)=0$ where
    \[f(s) := \phi_1(s,y)+\log x,\]
    using Newton's Method. 
  We require that 
		$\abs{\alpha'-\alpha} < \min\{0.0001,1/(\overline{u}\log x)\}$ 
  where $\overline{u} = \min\{\log x/\log y, y/\log y\}$.
    \item Output $HT(x,y,\alpha')$.
\end{enumerate}
The overall running time is
$${O\left(y\left[\frac{\log\log x}{\log y}
  +\frac{1}{\log \log y}\right]\right)}.$$
We have
\[f(s) = \phi_1(s,y)+\log x
    \mbox{ and } f'(s) = \phi_2(s,y)\]
So, our iteration function $g$ for Newton's Method is
\[g(s) := s - f'(s)/f(s).\]
\begin{itemize}
	\item This algorithm has a running time of $O(\pi(y)\log\log x)$:
	    $O(\log\log x)$ iterations of Newton's method to converge, with
		each iteration requiring a sum over the primes $\le y$
		to evaluate $\phi_1$.
    \item In practice, 5-6 iterations suffice 
	    for Newton's Method to converge.
\end{itemize}

Algorithm HT-fast \cite{Sorenson2000} estimates the functions
  $\zeta,\phi_1,\phi_2$ using the prime number theorem; 
  the Riemann Hypothesis is used to bound the error.
  The error is also controlled by evaluating an initial segment over
  the primes up to $O(\sqrt{y})$.
\begin{itemize}
    \item Set $z := \mbox{min}\{y,\mbox{max}\{1000,5\sqrt{y}\}\}$.
	We have $\phi_1(s,y) \approx -B(s,y,z)$ where
    $$
        B(s,y,z) = \sum_{p\leq z} \frac{\log p}{p^s-1} 
	+ \sum_{k=1}^{\left \lfloor{(\log y)/s)}\right \rfloor } 
		\frac{1}{1-ks} \left( y^{1-ks}-z^{1-ks}\right).
    $$
\end{itemize}
The functions $\zeta,\phi_2$ are similarly approximated.
This version of the algorithm is much faster, taking time proporional to
  $\sqrt{y}$, but gives estimates that, though still good, 
   are not quite as accurate as Algorithm HT in practice.

\subsection{Our New Algorithm HT$\alpha$}

Here are the steps, following the idea outlined at the beginning of
this section.
\begin{enumerate}
    \item Find all primes $p\leq y$
    \item Starting at $\alpha_0 := \log(1+y/(5\log x))/\log y$, 
	compute an approximation
	$\alpha_f$ to the solution $s=\alpha_1$ of $f_1(s)=0$ where
    \[f_1(s) := -B(s,y,z)+\log x.\]
    We must have 
	$\abs{\alpha_1-\alpha_f} < \min\{0.000001,1/\overline{u}\log x\}$.
	%	and $\abs{\log x - B(\alpha_f,y,z)}\leq 1$
    \item Using $\alpha_f$ as a starting point,
	    compute the approximation $\alpha$ to the solution 
		$s=\alpha_2$ of $f_2(s)=0$ where
    \[f_2(s) := \phi_1(s,y)+\log x,\]
    as before.
		We must have $\abs{\alpha_2-\alpha} < \mbox{min}\{0.0001,1/(\overline{u}\log x)\}$
    \item Output $HT(x,y,\alpha)$
\end{enumerate}
The running times for each of the steps of Algorithm HT$\alpha$ are as follows:
\begin{enumerate}
    \item  $O(y/\log\log y)$
    \item $O(\sqrt{y}\log\log x/\log y)$
    \item $O(y/\log y)$ per iteration
    \item $O(y/\log y)$
\end{enumerate}
Observations:
\begin{itemize}
    \item We can often treat step (1) as a preprocessing step.  In practice,
	    we often already have a list of primes available.
    \item In practice, we observed that step (3) 
	    will only run for 1-2 iterations.
  \item
	  In theory, one iteration suffices for step (3)
		if $(\log y)^2\gg\log x$, 
		as the accuracy
		guarantee for HT-fast matches that of HT to within
		a factor of $(1+O(\log y/\log x + 1/\log y))$, which
		matches $(1+O(1/\overline{u}))$, the relative error
		of Algorithm HT, in this case.
	In general, however, $O(\log\log x)$ iterations may be required
		if $y$ is extremely small compared to $x$, but in
		this case Algorithm HT is already fast.
  \item In no situation should Algorithm HT$\alpha$ be slower than
	  Algorithm HT.
	\item Algorithm HT-fast relies on the Riemann Hypothesis (RH) for 
	  correctness, but HT$\alpha$ only relies on
	the RH for running time.
\end{itemize}
Thus, our overall running time is reduced to $O(\pi(y))$,
%$\displaystyle{O\left(\frac{y}{\log y}\right)}$, 
under the assumption we have a list of primes available and $y$
is not extremely small compared to $x$.

For references on prime sieves, see \cite{AB2004,Helfgott2020,Sorenson2015}.

\subsection{Implementation Results}

Algorithm HT$\alpha$ has comparable error to Algorithm $HT$ 
and a faster running time.
\begin{figure}[ht]
\begin{center}
  \begin{tabular}{ccc}
	  \includegraphics[width=0.4\textwidth]{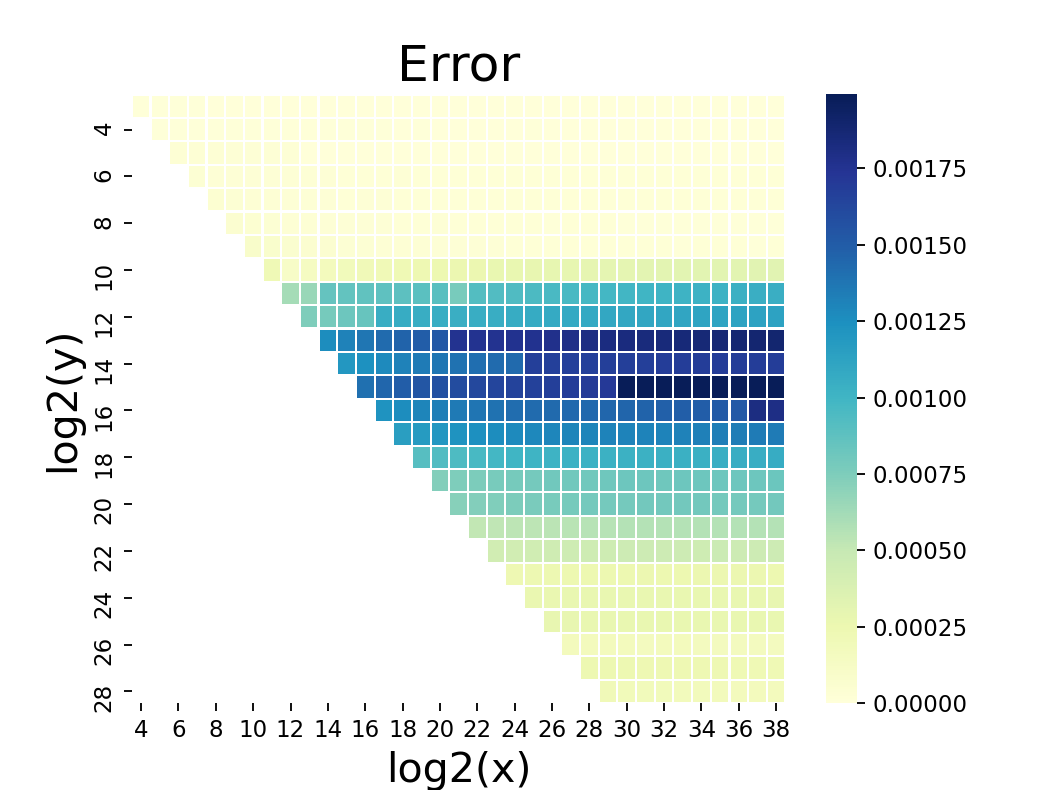} &
	  \hspace*{0.25in} &
	\includegraphics[width=0.4\textwidth]{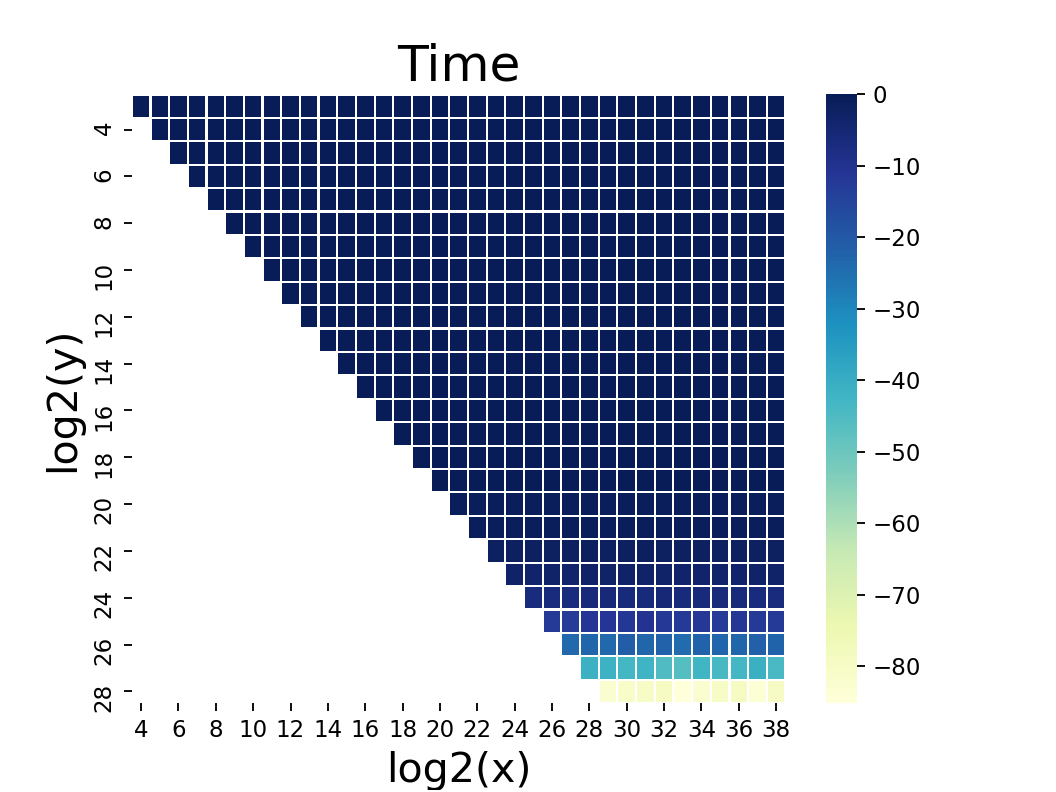} \\[5pt]
	  $\displaystyle{\abs{HT-HT\alpha}/\Psi(x,y)}$ &
	  \hspace*{0.25in} &
	  (HT time)-(HT$\alpha$ time) \\
  \end{tabular}
\end{center}
\caption{Performance Comparison: Algorithm HT versus HT$\alpha$
	\label{fig:alpha}}
\end{figure}

The following table is a comparison of error, 
time (in seconds) and Newton's Method iterations (Its.) per algorithm and
  by step in the case of HT$\alpha$.
Here $x = 2^{30}$:
\begin{center}
\begin{tabular}{|c||c|c|c||c|c|c|c|}
    \hline
    &\multicolumn{3}{|c||}{HT}& \multicolumn{4}{c|}{HT$\alpha$} \\ \hline
	$y$ & $HT/\Psi$ & Time & Its. & $HT\alpha/\Psi$ & Time & Its. (2) & Its. (3) \\ \hline
    \hline
     $2^{15}$ & 1.004 & 0.027 & 5 & 1.014 & 0.0006 & 6 & 1\\ \hline
    $2^{20}$ &1.031 & 0.765 & 6 & 1.034 & 0.002 & 6 & 1\\ \hline
    $2^{25}$ & 1.018 &  18.599 &  6 &  1.019 & 6.978 & 6 & 1\\ \hline
    \end{tabular}
\end{center}
The graphs in Figure \ref{fig:alpha} give a bit more data
  visually.

%%%%eval

\section{Estimating $\Psi(x,y)$\label{sec:eval}}

In this section, we give advice on the best way to compute
  values of $\Psi(x,y)$ for various ranges of $x$ and $y$ in practice.
\begin{figure}[ht]
\begin{center}
    \includegraphics[width=1.1\textwidth]{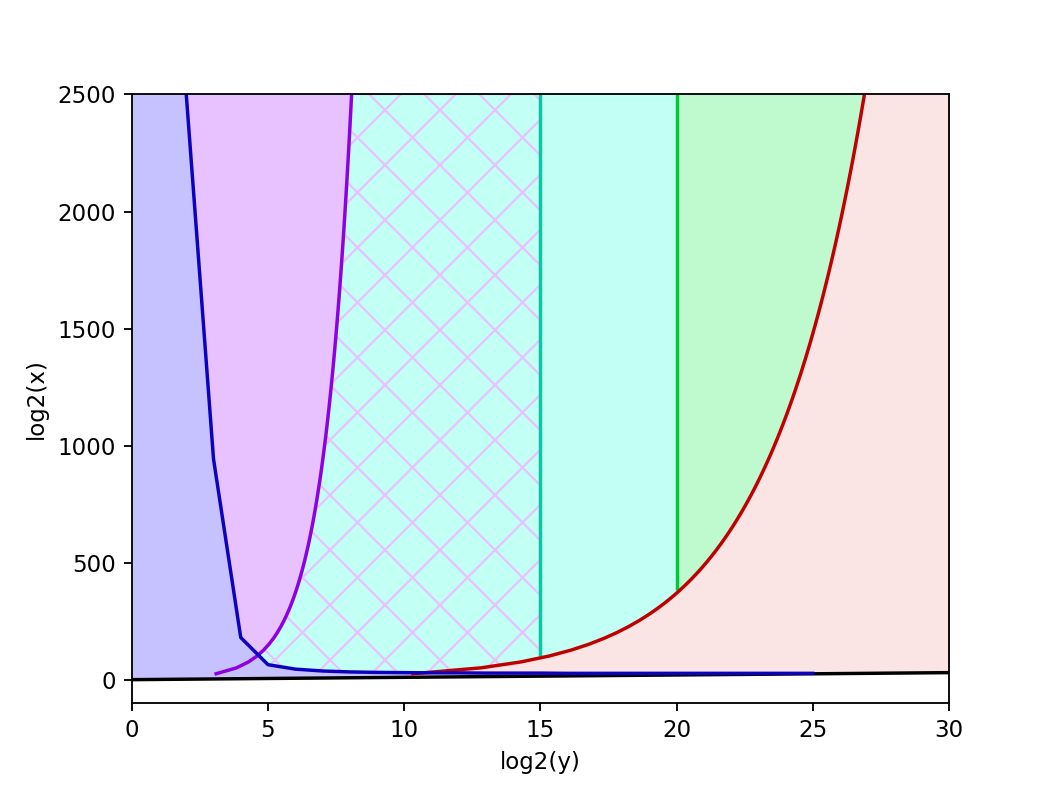}
\end{center}
	\caption{Algorithm Recommendations\label{fig:gp}}
\end{figure}
\begin{figure}[ht]
\begin{center}
    \includegraphics[width=1.1\textwidth]{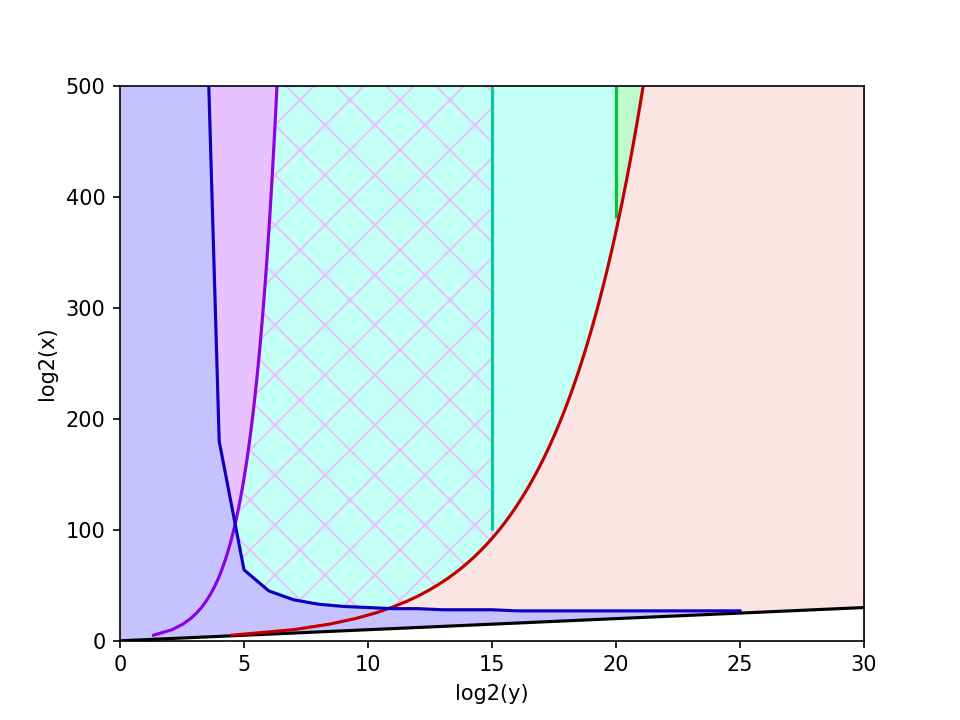}
\end{center}
	\caption{Algorithm Recommendations (small $x$)\label{fig:gpsmall}}
\end{figure}

We considered the following algorithms to estimate $\Psi(x,y)$,
  roughly ordered by how large $y$ is relative to $x$,
  starting with methods that work best for small $y$.
\begin{itemize}
\item \textbf{Buchstab's Identity}, stated below, directly implies
  a simple recursive algorithm that gives exact values
  of $\Psi(x,y)$.  However, with a running time roughly
  proportional to the value of $\Psi(x,y)$, it is
  only useful for very small inputs.
  \begin{eqnarray*}
    \Psi(x,y)&=&1+\sum_{p\leq y} \Psi(x/p,p)
  \end{eqnarray*}
  In practice, we supplement with the base case 
  $\Psi(x,2)=\lfloor \log_2 x \rfloor +1$ as well.

\item \textbf{Ennola's Second Theorem} was discussed in 
  detail in \S\ref{sec:ennola}.
  It is provably useful for when $y\le (\log x)^{3/4-\epsilon}$,
  but in practice we found it to be perfectly fine so long a its 
  running time is tolerable.
  Precomputation requires an $O(y^2/\log y)$ running time, which
  is quite high, but if this information is saved, then the time
  to compute specific $\Psi(x,y)$ values drops to a more reasonable
  $O(y\log y)$ time.
  It is highly accurate in practice.

\item \textbf{Algorithm HT} (or HT$\alpha$) was discussed above in
  \S\ref{sec:htalpha}.
  This method is provably accurate for $2\le y \le x$, but is a bit
  slow with a running time roughly linear in $y$.
  Its running time is similar to Ennola's method if precomputation is
  allowed, and much faster if not. 
  Also, less extra space is required for Algorithm HT.
  Ennola's method is a bit more accurate, but not provably so.

\item \textbf{Algorithm HT-fast} is the version of the previous
  algorithm where sums of primes are estimated using the Riemann
  Hypothesis to bound the error.
  It is a bit less accurate than Algorithm HT, but much faster with
  a running time proportional to $\sqrt{y}$ and the same wide range
  of applicability.

\item The \textbf{Dickman $\rho$} estimate gives 
    \[\Psi(x,y) \approx x\cdot \rho(u) + (1-\gamma)\frac{x}{\log x} \rho(u-1)\]
  where $u = \log x/\log y$ and $\rho(u)$ is the unique solution the
  the following equations:
  \begin{eqnarray*}
	  \rho(u)& =& 1  \quad (0\leq u \leq 1), \\
	  \rho'(u)& =& -\rho(u-1)/u  \quad (u\geq 1).
  \end{eqnarray*}
  Note that in the literature, one normally sees
  $$ \Psi(x,y) \sim x\cdot \rho(u), $$
  but we find adding the second term is worthwhile in improved
  accuracy.
  This estimate is valid when $y \geq L(x)$, with
    $L(x) = (\log x)^{2+\varepsilon}$
  assuming the Extended Riemann Hypothesis (ERH).
  Without the ERH, the lower limit on $y$ is much larger,
  $ \exp (\log_2 x)^{5/3+\epsilon}  $
  for $\epsilon>0$ \cite{Hildebrand86}.
  This method is very fast; with precomputaion of the $\rho$ function,
  evaluations of $\Psi(x,y)$ take constant time.
  $\rho(u)$ can be computed reasonably quickly using numeric integration;
		see \cite{vdLW69,HS97}.
  It is the least accurate of the methods presented here, but its
  very fast computing time make it desireable, especially for large $y$,
  where its accuracy is tolerable in practice.
\end{itemize}
In Figures \ref{fig:gp} and \ref{fig:gpsmall} are plots indicating 
approximately where we recommend one uses each method to compute $\Psi(x,y)$.
They are plotted using a logarithmic scale (base 2) in both $x$ and $y$.
The key to the graphs is as follows:
\begin{center}
\begin{tabular}{cl|cl}
	\includegraphics[width = 0.075\textwidth, height=6mm]{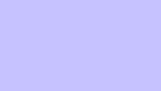}&  Use Buchstab's identity  &
     \includegraphics[width = 0.025\textwidth, height=6mm]{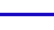}& 
	Buchstab takes 1.5 seconds\\
	\includegraphics[width = 0.075\textwidth, height=6mm]{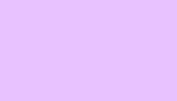}
	&  Use Ennola's Second Theorem &
	\includegraphics[width = 0.025\textwidth, height=6mm]{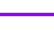}
	& $y=(\log x)^{3/4}$ \\
	\includegraphics[width = 0.075\textwidth, height=6mm]{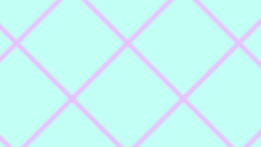}
	&  Use HT (Ennola's is accurate too) &
     \includegraphics[width = 0.025\textwidth, height=6mm]{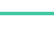} 
	& Ennola max (gets slow)  \\
	\includegraphics[width = 0.075\textwidth, height=6mm]{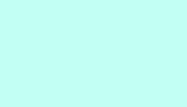}
	& Use HT (or HT$\alpha$) &
     \includegraphics[width = 0.025\textwidth, height=6mm]{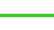}
	&  Switch from HT to HT-Fast \\
	\includegraphics[width = 0.075\textwidth, height=6mm]{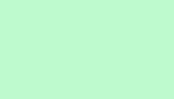}
	& Use HT-Fast &
     \includegraphics[width = 0.025\textwidth, height=6mm]{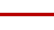}
	& $y=L(x) = (\log x)^{2.5}$ \\
	\includegraphics[width = 0.075\textwidth, height=6mm]{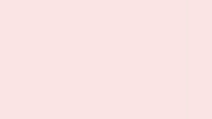}
	& Use $x\cdot \rho(u) + (1-\gamma) \frac{x}{\log x} \rho(u-1)$ &
	\includegraphics[width = 0.025\textwidth, height=6mm]{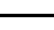}
	&  $y = x$ \\
\end{tabular}

\end{center}
\subsection*{Methodology and Notes}
\begin{itemize}
    \item We implemented all the algorithms in C++, 
	    we used the standard GNU g++ compiler,
	    and the code was run on
	    a Linux server using standard Intel hardware.
    \item In \cite{HS97} it was shown that Algorithm HT is extremely accurate
	    in practice, and so we used either Buchstab's algorithm, or
	Algorithm HT when Buchstab was too slow, as a baseline for accuracy.
    \item We collected data on the accuracy and speed of all the algorithms
	   over the range of $x,y$ values shown in the graphs, except
		for some algorithms that got too slow.
    \item After the data was collected, we determined the $x,y$ ranges
	    where each algorithm was reasonably accurate in practice.
    \item We varied the value of $\epsilon$ in the ERH cutoff
	    $y=L(x)=(\log x)^{2+\epsilon}$ for the Dickman $\rho$ method
		to see what worked best in practice.
\end{itemize}

%%%%wrapup

\section{Concluding Remarks\label{sec:wrapup}}
\begin{itemize}
\item
  There are algorithms that give explicit upper and lower bounds on
  $\Psi(x,y)$; see \cite{Bernstein98,Bernstein2002,LP2018,PS06}.
  Such methods tend to be noticably slower than the methods used here,
    which is why we did not consider them.
  That said, they do have their purposes.
\item
  We currently have no theory as to why the algorithm based on Ennola's
  second theorem has, in practice, applicability on a range as wide
  as that of Algorithm HT.  It may be that a closer examination
  of Ennola's proof may yield a way to improve the range.
\item
  Is it possible to speed up the algorithm from \S\ref{sec:ennola} 
  using the prime
  number theorem to estimate the sums and products over primes,
  perhaps bounding the error using the ERH as was done
  in \cite{PS06,Sorenson2000}?
\item
  In \cite{BS13} it was shown how to use LMO summation to
  improve the running time of Algorithm HT to $y^{2/3+\epsilon}$.
  It stands to reason that HT-fast can be done in time $y^{1/3+\epsilon}$
  as well.
  As far as the authors are aware, this has not yet been implemented.
\end{itemize}

\section*{Acknowledgements}

The first author was supported in part by the 
Butler Summer Institute (Summer 2021)
and the Mathematics Research Camp (August 2021) at Butler University.

This work was presented by the first author at
the Young Mathematicians Conference at The Ohio State University, August 2021,
and at the Undergraduate Research Conference at Butler University,
April 2022.

%\bibliography{all}

\begin{thebibliography}{10}

\bibitem{AB2004}
A.~O.~L. Atkin and D.~J. Bernstein.
\newblock Prime sieves using binary quadratic forms.
\newblock {\em Mathematics of Computation}, 73:1023--1030, 2004.

\bibitem{BS}
Eric Bach and Jeffrey~O. Shallit.
\newblock {\em Algorithmic Number Theory}, volume~1.
\newblock MIT Press, 1996.

\bibitem{BS2020}
Eric Bach and Jonathan Sorenson.
\newblock An algorithm to generate random factored smooth integers, 2020.
\newblock Available on arxiv.org at https://arxiv.org/abs/2006.07445.

\bibitem{BS13}
Eric Bach and Jonathan~P. Sorenson.
\newblock Approximately counting semismooth integers.
\newblock In {\em Proceedings of the 38th International symposium on symbolic
  and algebraic computation}, ISSAC '13, pages 23--30, New York, NY, USA, 2013.
  ACM.

\bibitem{Bernstein95-2}
Daniel~J. Bernstein.
\newblock Enumerating and counting smooth integers.
\newblock Chapter 2, PhD Thesis, University of California at Berkeley, May
  1995.

\bibitem{Bernstein98}
Daniel~J. Bernstein.
\newblock Bounding smooth integers.
\newblock In J.~P. Buhler, editor, {\em Third International Algorithmic Number
  Theory Symposium}, pages 128--130, Portland, Oregon, June 1998. Springer.
\newblock LNCS 1423.

\bibitem{Bernstein2002}
Daniel~J. Bernstein.
\newblock Arbitrarily tight bounds on the distribution of smooth integers.
\newblock In Bennett, Berndt, Boston, Diamond, Hildebrand, and Philipp,
  editors, {\em Proceedings of the Millennial Conference on Number Theory},
  volume~1, pages 49--66. A. K. Peters, 2002.

\bibitem{Couveignes2006}
J.~Couveignes.
\newblock Hard homogeneous spaces.
\newblock 2006.

\bibitem{CP}
R.~Crandall and C.~Pomerance.
\newblock {\em Prime Numbers, a Computational Perspective}.
\newblock Springer, 2001.

\bibitem{Ennola69}
V.~Ennola.
\newblock On numbers with small prime divisors.
\newblock {\em Ann. Acad. Sci. Fenn. Ser. A I}, 440, 1969.
\newblock 16pp.

\bibitem{Granville2008}
Andrew Granville.
\newblock Smooth numbers: computational number theory and beyond.
\newblock In {\em Algorithmic number theory: lattices, number fields, curves
  and cryptography}, volume~44 of {\em Math. Sci. Res. Inst. Publ.}, pages
  267--323. Cambridge Univ. Press, Cambridge, 2008.

\bibitem{Helfgott2020}
Harald~Andr\'{e}s Helfgott.
\newblock An improved sieve of {Eratosthenes}.
\newblock {\em Mathematics of Computation}, 89:333--350, 2020.

\bibitem{Hildebrand86}
A.~Hildebrand.
\newblock On the number of positive integers $\le x$ and free of prime factors
  $> y$.
\newblock {\em Journal of Number Theory}, 22:289--307, 1986.

\bibitem{HT86}
A.~Hildebrand and G.~Tenenbaum.
\newblock On integers free of large prime factors.
\newblock {\em Trans. AMS}, 296(1):265--290, 1986.

\bibitem{HT93}
A.~Hildebrand and G.~Tenenbaum.
\newblock Integers without large prime factors.
\newblock {\em Journal de Th\'eorie des Nombres de Bordeaux}, 5:411--484, 1993.

\bibitem{HS97}
Simon Hunter and Jonathan~P. Sorenson.
\newblock Approximating the number of integers free of large prime factors.
\newblock {\em Mathematics of Computation}, 66(220):1729--1741, 1997.

\bibitem{KTP76}
D.~E. Knuth and L.~Trabb~Pardo.
\newblock Analysis of a simple factorization algorithm.
\newblock {\em Theoretical Computer Science}, 3:321--348, 1976.

\bibitem{LP2018}
J.~D. Lichtman and C.~Pomerance.
\newblock Explicit estimates for the distribution of number free of large prime
  factors.
\newblock {\em Journal of Number Theory}, 183:1--23, 2018.

\bibitem{MOV}
A.~J. Menezes, P.~C. van Oorschot, and S.~A. Vanstone.
\newblock {\em Handbook of Applied Cryptography}.
\newblock CRC Press, Boca Raton, 1997.

\bibitem{PS06}
Scott Parsell and Jonathan~P. Sorenson.
\newblock Fast bounds on the distribution of smooth numbers.
\newblock In Florian Hess, Sebastian Pauli, and Michael Pohst, editors, {\em
  Proceedings of the 7th International Symposium on Algorithmic Number Theory
  (ANTS-VII)}, pages 168--181, Berlin, Germany, July 2006. Springer.
\newblock LNCS 4076, ISBN 3-540-36075-1.

\bibitem{RS2006}
A.~Rostovtsev and A.~Stoblunov.
\newblock Public-key cryptosystem based on isogenies.
\newblock 2006.

\bibitem{Sorenson2000}
Jonathan~P. Sorenson.
\newblock A fast algorithm for approximately counting smooth numbers.
\newblock In W.~Bosma, editor, {\em Proceedings of the Fourth International
  Algorithmic Number Theory Symposium (ANTS IV)}, pages 539--549, Leiden, The
  Netherlands, 2000.
\newblock LNCS 1838.

\bibitem{Sorenson2015}
Jonathan~P. Sorenson.
\newblock Two compact incremental prime sieves.
\newblock {\em LMS Journal of Computation and Mathematics}, 18(1):675--683,
  2015.

\bibitem{Suzuki2004}
K.~Suzuki.
\newblock An estimate for the number of integers without large prime factors.
\newblock {\em Mathematics of Computation}, 73:1013--1022, 2004.
\newblock MR 2031422 (2005a:11142).

\bibitem{Suzuki2006}
K.~Suzuki.
\newblock Approximating the number of integers without large prime factors.
\newblock {\em Mathematics of Computation}, 75:1015--1024, 2006.

\bibitem{Tenenbaum}
G\'erald. Tenenbaum.
\newblock {\em Introduction to Analytic and Probabilistic Number Theory},
  volume~46 of {\em Cambridge Studies in Advanced Mathematics}.
\newblock Cambridge University Press, english edition, 1995.

\bibitem{vdLW69}
J.~van~de Lune and E.~Wattel.
\newblock On the numerical solution of a differential-difference equation
  arising in analytic number theory.
\newblock {\em Mathematics of Computation}, 23:417--421, 1969.

\end{thebibliography}

\end{document}